\pdfoutput=1
\RequirePackage{ifpdf}
\ifpdf 
\documentclass[pdftex]{sigma}
\else
\documentclass{sigma}
\fi

\usepackage{wrapfig,case,txmac}

\numberwithin{equation}{section}
\newtheorem{theo}{Theorem}[section]
\newtheorem{lem}[theo]{Lemma}
\newtheorem{propos}[theo]{Proposition}
{\theoremstyle{definition}
\newtheorem{defn}[theo]{Definition}
\newtheorem{rem}[theo]{Remark}
}

\DeclareMathOperator{\tr}{tr}
\DeclareMathOperator{\spaan}{span}

\DeclareMathOperator{\spn}{span}

\begin{document}

\def\|{\mathrel{\kern1.5pt\Vert\kern1.5pt}}
\def\br#1{\left\langle#1\right\rangle}
\def\iv#1{\left\lfloor#1\right\rfloor}

\allowdisplaybreaks

\renewcommand{\PaperNumber}{105}

\FirstPageHeading

\ShortArticleName{Everywhere Equivalent 3-Braids}

\ArticleName{Everywhere Equivalent 3-Braids}

\Author{Alexander STOIMENOW}

\AuthorNameForHeading{A.~Stoimenow}

\Address{Gwangju Institute of Science and Technology, School of General Studies, GIST College,\\
123 Cheomdan-gwagiro, Gwangju 500-712, Korea}
\Email{\href{mailto:stoimeno@stoimenov.net}{stoimeno@stoimenov.net}}
\URLaddress{\url{http://stoimenov.net/stoimeno/homepage/}}

\ArticleDates{Received July 08, 2014, in f\/inal form November 04, 2014; Published online November 16, 2014}

\Abstract{A knot (or link) diagram is said to be everywhere equivalent if all the diagrams obtained by switching one
crossing represent the same knot (or link).
We classify such diagrams of a~closed 3-braid.}


\Classification{57M25; 20F36; 20E45; 20C08}

\section{Introduction}

How does a~diagram~$D$ of a~knot (or link)~$L$ look, which has the following property: all diag\-rams~$D'$ obtained~by
changing exactly one crossing in~$D$ represent the same knot (or link) $L'$ (which we allow to be dif\/ferent from~$L$)?
This suggestive question was possibly f\/irst proposed in this form by K.~Taniyama, who called such diagrams
\em{everywhere equivalent} (see Def\/inition~\ref{Md} in Section~\ref{S@}).
Through the connection between links and braid groups (see, e.g.,~\cite{Jones}), this question turns out related to the
following group-theoretic question: given a~group in a~set of conjugate generators, which words in these generators have
the property that reversing each individual letter gives a~set of conjugate elements?

When by the crossing-changed diagrams $D'$ display the unknot, Taniyama's question was studied previously
in~\cite{nikos}, where~$D$ was called \em{everywhere trivial}.
There some ef\/forts were made to identify such diagrams, mostly by computationally verifying a~number of low-crossing
cases.
The upshot was that, while there is a~(hard to be described) abundance of diagrams for~$D$ unknotted, only 6 simple
diagrams seem to occur when~$D$ is not; see~\eqref{eeq}.

Motivated by Taniyama's more general concept, we made in~\cite{evrdiff} a~further study of such phenomena.
We conjectured, as an extension of the everywhere trivial case, a~general description of everywhere equivalent diagrams
for a~knot, and proved some cases of low genus diagrams.
We also proposed some graph-theoretic constructions of everywhere equivalent diagrams for links.

In this paper we give the answer for 3-braids.
A~3-braid diagram corresponds to a~particular braid \em{word} in Artin's generators.
This word can be regarded up to inversion, cyclic permutations, the interchanges $\sigma_i\leftrightarrow\sigma_i^{-1}$
(ref\/lection) and $\sigma_1\leftrightarrow\sigma_2$ (f\/lip).
However, beyond this it will be of importance to distinguish between braids and their words, i.e., how a~given braid is written.

\begin{theo}
\label{th1}
A~$3$-braid word~$\beta$ gives an everywhere equivalent diagram if and only if it is in one of the following four families
{\rm (}up to the above equivalence$)$:
\def\labelenumi{{\rm \theenumi)}}
\begin{enumerate}\itemsep=0pt
\item $\big(\sigma_1\sigma_2^{-1}\big)^k$ or $\big(\sigma_1\sigma_2\sigma_1^{-1}\sigma_2^{-1}\big)^{k}$ for $k=1,2$,
and $\sigma_1\sigma_2^{-1}\sigma_1^{-2}\sigma_2^2$ {\rm (non-positive} case; refers to all three possibilities$)$,
\label{type1}
\item {\rm any} positive $($or negative, including the trivial$)$ word representing a~central element
$\Delta^{2k}$, $k\in{\mathbb Z}$ {\rm (central} case$)$,
\label{type2}
\item the words $\big(\sigma_1^l\sigma_2^l\big)^k$ for $k,l\ge 1$ {\rm (symmetric} case$)$, and
\label{type3}
\item the words $\sigma_1^k$ for $k>0$ {\rm (split} case$)$.
\label{type4}
\end{enumerate}
\end{theo}

Using the exponent sum~\eqref{xs}, one easily sees that the answer to the related group-theoretic question for the
3-braid group consists of the last three families in the theorem.
This outcome, and even more so its derivation, have turned out more complicated than expected.

The f\/irst family mainly comes, except for $\sigma_1\sigma_2\sigma_1^{-1} \sigma_2^{-1}$ and the trivial cases, from the
diagrams of~\eqref{eeq} (under the exclusion of the 5 crossing one, which is not a~3-braid diagram).
The last two families are also quite suggestive, and given (in more general form for arbitrary diagrams)
in~\cite{evrdiff}.
However, we initially entirely overlooked the second family of diagrams.
We did not explicitly ask in~\cite{evrdiff} whether our link diagram constructions are exhaustive, but we certainly had
them in mind when approaching Theorem~\ref{th1}.
These previous examples come in some way from totally symmetric (edge transitive) planar graphs.
However, there is little symmetry in general here.
One can easily construct examples of central (element) words lacking any symmetry.
(One can also see that every positive word in the 3-braid group can be realized as a~subword of a~positive word
representing a~central element.)

The second family does not yield (and thus answer negatively our original question for) knots, and it does not lead out
of the positive case (whose special role was well recognized in~\cite{evrdiff}).
Still we take it as a~caution that everywhere equivalence phenomena, although sporadic, may occur in far less
predictable ways than believed.

Our proof is almost entirely algebraic, and will consist in using the Jones polynomial to distinguish the links of
various $D'$ except in the desired cases.
Mostly we will appeal to the description of the Jones polynomial in terms of the Burau representation, but at certain
points it will be important to use information coming from the skein relation and the Kauf\/fman bracket.

The proof occupies Sections~\ref{S3} and~\ref{S4}, divided by whether the braid word is positive or not.
Both parts require somewhat dif\/ferent treatment.
We will see that for 3-braids the non-positive case quickly recurs to the everywhere trivial one.

A f\/inal observation (Proposition~\ref{p2}) addresses the lack of interest of the situation opposite to everywhere
equivalence: when all crossing-switched versions of a~diagram are to represent dif\/ferent links.
(This property was called \em{everywhere different}, and some constructions for such knot diagrams were given
in~\cite{evrdiff}.)

In a~parallel paper~\cite{e12cp} we observe how to solve the classif\/ication of (orientedly) everywhere equivalent
diagrams in another case, this of two components.

\section{Preliminaries}

It seems useful to collect various preliminaries, which will be used at dif\/ferent places later in the paper.

\subsection{Link diagrams and Jones polynomial}

All link diagrams are considered oriented, even if orientation is sometimes ignored.
We also assume here that we actually regard the plane in which a~link diagram lives as $S^2$, that is, we consider as
equivalent diagrams which dif\/fer by the choice of the point at inf\/inity.

The \em{Jones polynomial}~$V$ can be def\/ined as the polynomial taking the value $1$ on the unknot, and satisfying the
\em{skein relation}
\begin{gather}
\label{1}
t^{-1}V\Bigl( \text{\raisebox{-6pt}[0pt][0pt]{\includegraphics{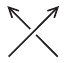}}} \Bigr)-t V\Bigl(\text{\raisebox{-6pt}[0pt][0pt]{\includegraphics{Stoimenow-Pic1}}}  \Bigr) = \big(t^{1/2}-t^{-1/2}\big)V\Bigl(\text{\raisebox{-5pt}[0pt][0pt]{\includegraphics{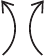}}} \Bigr).
\end{gather}
In each triple as in~\eqref{1} the link diagrams are understood to be identical except at the designated spot.
The fragments are said to depict a~\em{positive crossing}, \em{negative crossing}, and \em{smoothed out} crossing.
The \em{skein smoothing} is thus the replacement of a~crossing by the third fragment.
The \em{writhe} $w(D)$ of a~link diagram~$D$ is the sum of the signs of all crossings of~$D$.
If all crossings of~$D$ are positive, then~$D$ is called positive.

\begin{figure}[htb]\centering
\includegraphics{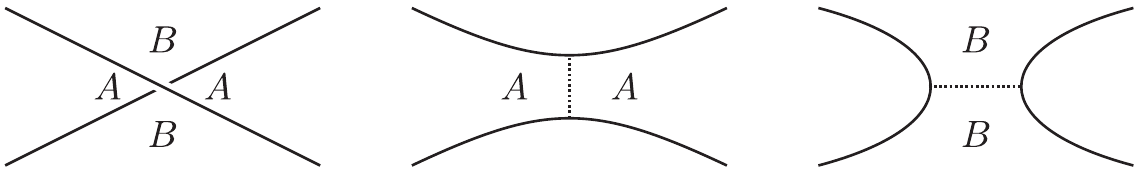}
\caption{The $A$- and $B$-corners of a~crossing, and its both splittings.
The corner~$A$ (resp.~$B$) is the one passed by the overcrossing strand when rotated counterclockwise (respectively
clockwise) towards the undercrossing strand.
A~type~$A$ (resp.~$B$) splitting is obtained by connecting the~$A$ (resp.~$B$) corners of the crossing.
The dashed line indicates the trace of the crossing after the split.}\label{figsplit}
\end{figure}

It is useful to recall here the alternative description of~$V$ via Kauf\/fman's state model~\cite{Kauffman}.
A~\em{state} is a~choice of \em{splicings} (or \em{splittings}) of type~$A$ or~$B$ (see Fig.~\ref{figsplit}) for any
single crossing of a~link diagram~$D$.
We call the \em{$A$-state} the state in which all crossings are~$A$-spliced, and the \em{$B$-state} $B(D)$ is def\/ined
analogously.

When for a~state~$S$ all splicings are performed, we obtain a~\em{splicing diagram}, which consists of a~collection of
(disjoint) \em{loops} in the plane (solid lines) together with \em{$($crossing$)$ traces} (dashed lines).
We call a~loop \em{separating} if both its interior and exterior contain other loops (regardless of what traces).
We will for convenience identify below a~state~$S$ with its splicing diagram for f\/ixed~$D$.
We will thus talk of the loops and traces of a~state.

Recall, that the \em{Kauffman bracket} $\br{D}$~\cite{Kauffman} of a~link diagram~$D$ is a~Laurent polynomial in
a~variable~$A$, obtained by a~sum over all states~$S$ of~$D$:
\begin{gather}
\label{eq_12}
\br{D} =\sum\limits_S A^{\#A(S)-\#B(S)} \big({-}A^2-A^{-2}\big)^{|S|-1}.
\end{gather}
Here $\#A(S)$ and $\#B(S)$ denote the number of type~$A$ (respectively, type~$B$) splittings and $|S|$ the number of (solid
line) loops in the splicing diagram of~$S$.
The formula~\eqref{eq_12} results from applying the f\/irst of the \em{bracket relations}
\begin{gather*}
\Big\langle\text{\raisebox{-6pt}[0pt][0pt]{\includegraphics{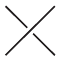}}}\Big\rangle
=
A^{-1}\Big\langle\text{\raisebox{-4pt}[0pt][0pt]{\includegraphics{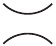}}}\Big\rangle
+A\Big\langle\text{\raisebox{-4.5pt}[0pt][0pt]{\includegraphics{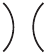}}}\Big\rangle ,
\qquad
\Big\langle\text{\raisebox{-5pt}[0pt][0pt]{\includegraphics{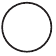}}} \cup X\Big\rangle
=
\big({-}A^2-A^{-2}\big)\br{X},
\end{gather*}
to each crossing of~$D$ (here traces are ignored), and then deleting (except one) loops using the second relation, at
the cost of a~factor $-A^2-A^{-2}$ per deleted loop.
(The normalization is thus here that the diagram of one circle with no crossings has unit bracket.)

The Jones polynomial of a~link~$L$ can be determined from the Kauf\/fman bracket of some diagram~$D$ of~$L$~by
\begin{gather}
\label{conv}
V_L(t) =\big({-}t^{-3/4}\big)^{-w(D)}\br{D} \raisebox{-0.6em}{$
\Big|_{A=t^{-1/4}}$},
\end{gather}
with $w(D)$ being the writhe of~$D$.
This is another way, dif\/ferent from~\eqref{1}, to specify the Jones polynomial.

It is well-known that $V\in {\mathbb Z}[t^{\pm 1}]$ (i.e., only integral powers occur) for odd number of link components
(in particular, for knots), while $V\in t^{1/2}\cdot {\mathbb Z}[t^{\pm 1}]$ (i.e., only half-integral powers occur) for
even number of components.

For $V\in{\mathbb Z}[t^{1/2},t^{-1/2}]$, the \em{minimal} or \em{maximal degree} $\min\deg V$ or $\max\deg V$ is the
minimal resp.\ maximal exponent of~$t$ with non-zero coef\/f\/icient in~$V$.
Let $\spn V=\max\deg V-\min\deg V$.

\subsection{Semiadequacy and adequacy}

Let~$S$ be the~$A$-state of a~diagram~$D$ and $S'$ a~state of~$D$ with exactly one~$B$-splicing.
If $|S|>|S'|$ for all such $S'$, we say that~$D$ is \em{$A$-adequate}.
Similarly one def\/ines a~$B$-adequate diagram~$D$ (see~\cite{LickThis}).
Then we set a~diagram to be
\begin{gather}
\text{adequate}  =  \text{$A$-adequate}\enspace \text{\em and}\enspace\text{$B$-adequate},\nonumber\\
\text{semiadequate}  =  \text{$A$-adequate}\enspace\text{\em or}\enspace\text{$B$-adequate},\label{aqdf}  \\
\text{inadequate}  =  \text{\em neither}\enspace \text{$A$-adequate}\enspace\text{\em nor}\enspace\text{$B$-adequate}.\nonumber
\end{gather}
(Note that inadequate is a~stronger condition than not to be adequate.)

A link is called~$A$ (or~$B$)-adequate, if it has an~$A$ (or~$B$)-adequate diagram.
It is \em{semiadequate} if it is~$A$- or~$B$-adequate, and \em{inadequate}, if it is not semiadequate, that is,
neither~$A$- nor~$B$-adequate.
A~link is \em{adequate} if it has an adequate diagram.
This property is stronger than being both~$A$- and~$B$-adequate, since a~link might have diagrams that enjoy either
properties, but none that does so \em{simultaneously}.
The Perko knot $10_{161}$ in~\cite[Appendix]{Rolfsen} is an example.
This property of the Perko knot follows from work of Thistlethwaite~\cite{Thistle}, and is explained, e.g., in
Cromwell's book~\cite[p.~234]{Cromwell}, or (along with further examples) in~\cite{adeq}.

A link diagram~$D$ is said to be \em{split}, if its planar image (forgetting crossing information) is a~disconnected set.
A~\em{region} of~$D$ is a~connected component of the complement of this planar image.
At every crossing of~$D$ four regions meet; if two coincide, we call the crossing \em{nugatory}.
A~diagram with no nugatory crossings is \em{reduced}.

It is easily observed (as in~\cite{LickThis}) that a~reduced alternating diagram (and hence an alternating link) is
adequate, and that the~$A$- and~$B$-state of an alternating non-split diagram have no separating loops.

\looseness=-1
The maximal degree of~$A$ of the summands in~\eqref{eq_12} is realized by the~$A$-state~$S$.
However,~its contribution may be cancelled by this of some other state.
One situation when this does \em{not} happen is when~$D$ is~$A$-adequate.
Then the~$A$-state gives in~\eqref{eq_12} the unique contribution to the maximal degree in~$A$, and, via~\eqref{conv},
the minimal degree of~$V$.
We call this the \em{extreme~$A$-term}.
Thus, for~$A$-adequate diagrams, $\min\deg V$ can be read of\/f from the~$A$-state directly, and the coef\/f\/icient is $\pm 1$.
For not~$A$-adequate diagrams, the situation is a~bit more subtle and studied in~\cite{BaeMor}.

We will use the following important special case of that study.
When~$D$ is not~$A$-adequate, the~$A$-state has a~trace connecting a~loop to itself, which we call \em{self-trace}.
For a~given loop, we call a~pair of self-traces ending on it \em{intertwined} if they have the mutual position
\begin{gather}
\label{lpt}
\begin{split}
& \includegraphics{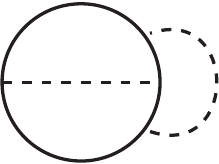}
\end{split}
\end{gather}
A~self-trace is \em{isolated} if it is intertwined with no other self-trace.
Bae and Morton show (among others) that if in the~$A$-state a~self-trace is isolated, the contribution in the
extreme~$A$-term of~$V$ is zero.

Similar remarks apply on~$B$-adequate diagrams and $\max\deg V$, and then on adequate diagrams and $\spn V$.

\subsection{Braid groups and words}\label{22}

The \em{$n$-string braid}, or shortly \em{$n$-braid, group} $B_n$ is considered generated by the Artin \em{standard
generators} $\sigma_i$ for $i=1, \dots,n-1$.
An Artin generator $\sigma_i$, respectively, its inverse $\sigma_i^{-1}$ will be called a~\em{positive}, respectively,
\em{negative} letter, and~$i$ the \em{index} of this letter.
We will almost entirely focus on $n=3$.

The Artin generators are subject to commutativity relations (for $n\ge 4$; the bracket denotes the commutator below) and
braid relations, which give $B_n$ the presentation
\begin{gather*}
B_n = \br{\sigma_1,\dots,\sigma_{n-1}\enspace\Bigg|
\begin{array}{ll}
[\sigma_i,\sigma_j]=1, & |i-j|>1,
\\
\sigma_{j}\sigma_i\sigma_{j}= \sigma_i\sigma_{j}\sigma_i, & |i-j|=1
\end{array}
}.
\end{gather*}

For the sake of legibility, we will commonly use a~bracket notation for braid words.
The meaning of this notation is the word obtained by replacing in the content of the brackets every integer $\pm i$ for
$i>0$, not occurring as an exponent, by $\sigma_i^{\pm 1}$, and removing the enclosing brackets.
Thus, e.g.,
\begin{gather*}
\big[1(2-1)^4-2-1\big] =\sigma_1\big(\sigma_2\sigma_1^{-1}\big)^4\sigma_2^{-1}\sigma_1^{-1}.
\end{gather*}
Although negative exponents will not be used much here, let us f\/ix for clarity that for a~letter they are understood as
the inverse letter, and for a~longer subword as the inverse letters written in reverse order.
Thus
\begin{gather*}
[-1-12]=\big[1^{-2}2\big]=\sigma_1^{-2}\sigma_2
\end{gather*}
and
\begin{gather*}
\big[(-12)^{-2}1\big]=\big[(-21)^21\big]=[-21-211]=\big[{-}21-21^2\big]=\sigma_2^{-1}\sigma_1 \sigma_2^{-1}\sigma_1^2.
\end{gather*}
Occasionally we will insert into the bracket notation vertical bars `$|$'.
They have no inf\/luence on the value of the expression, but we use them to highlight special subwords.

A word which does not contain a~letter followed or preceded by its inverse is called \em{reduced}.

\looseness=-1
In a~reduced braid word, a~maximal subword $\sigma_i^{\pm k}$ for $k>0$, i.e., one followed and preceded by letters of
dif\/ferent index, is called a~\em{syllable}.
The number~$i$ is called \em{index} of the syllable, and the number $\pm k$ its \em{exponent}, which is composed of its
\em{sign} `$\pm$' and its \em{length} $k>0$.
According to the sign, a~syllable is \em{positive} or \em{negative}.
A~syllable of length $1$ (of either sign) will be called \em{trivial}.

Obviously every reduced braid word decomposes into syllables in a~unique way:
\begin{gather*}
\beta  =\prod\limits_{i=1}^n\sigma_i^{p_i},
\end{gather*}
with $p_i\in{\mathbb Z}\setminus\{0\}$.
The sequence $(p_1,\dots,p_n)$ will be called \em{exponent vector}.
Thus an entry `$\pm 1$' in the exponent vector corresponds to a~trivial syllable.
A~word is \em{positive}, if all entries in its exponent vector are positive (i.e., it has no negative syllable).

Often braid words will be considered up to cyclic permutations.
In this case so will be done with the exponent vector.
The length of the exponent vector considered up to cyclic permutations will be called \em{weight} $\omega(\beta)$ of~$\beta$.
If $\beta \in B_3$, then $\omega(\beta)$ is even, since indices~$1$ and~$2$ can only interchange, except when
$\omega(\beta)=1$ (and~$\beta$ is a~single syllable).

The quantity $\sum|p_i|$ is the \em{length} of the word~$\beta$ (and is, of course, dif\/ferent from its weight, unless all
syllables are trivial).
The length-zero word will be called \em{trivial word}.
The quantity
\begin{gather}
\label{xs}
[\beta] =\sum\limits_{i=1}^np_i
\end{gather}
is called \em{exponent sum} of~$\beta$, and is invariant of the braid (i.e., equal for dif\/ferent words of the same
braid), and in fact its conjugacy class.

The \em{half-twist element} $\Delta\in B_n$ is given~by
\begin{gather*}
\Delta=(\sigma_1\sigma_2\cdots\sigma_{n-1})(\sigma_1\cdots\sigma_{n-2})\cdots
(\sigma_1\sigma_2)\sigma_1,
\end{gather*}
and its square
\begin{gather*}
\Delta^2=(\sigma_1\sigma_2\cdots\sigma_{n-1})^n
\end{gather*}
is the generator of the \em{center} of $B_n$.

We will need mostly the group $B_3$, where~$\Delta$ has the two positive word representations $[121]$ and $[212]$.
We will use the notation $\overline{\mathord{\cdot}}$ for the involution $\sigma_1\leftrightarrow\sigma_2$ of $B_3$
induced by conjugacy with~$\Delta$.

\subsection{Braids and links}

There is a~well-known graphical representation for braids:
\begin{gather*}
\includegraphics{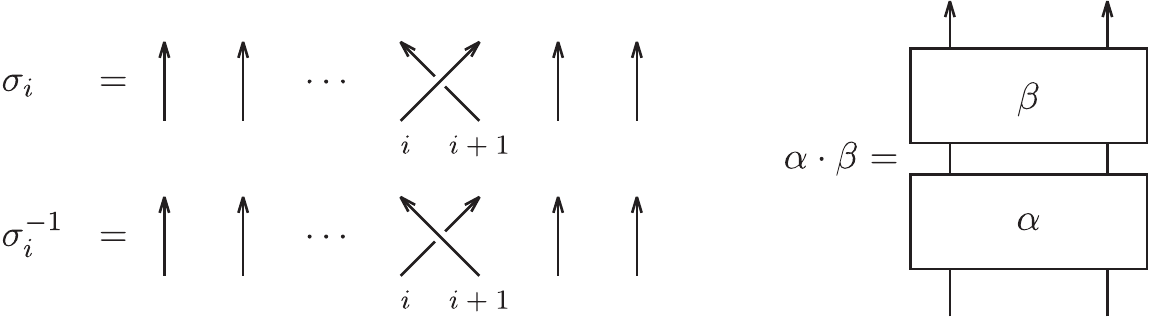}
\end{gather*}
Thus, a~(positive/negative) letter of a~braid word gives a~(positive/negative) crossing, and smoothing a~crossing
corresponds to deleting a~letter.
In this sense we will feel free to speak of (switching) crossings and smoothings of a~braid (word).

For braid(word)s~$\beta$ there is a~\em{closure} operation $\hat\beta $:
\begin{gather*}
\includegraphics{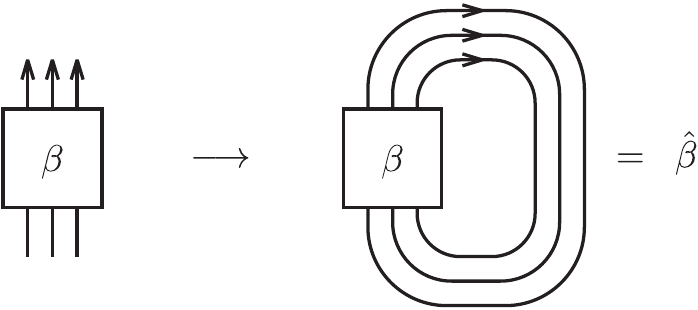}
\end{gather*}
In this way, a~\em{braid closes to a~knot or link} and a~\em{braid word} (which can be also regarded here up to cyclic
permutations) \em{closes to a~knot or link diagram}.
Here the separation between the two levels of correspondence must be kept in mind.

Thus~$\beta$ is a~positive word if and only if $\hat\beta $ is a~positive link diagram.
For a~further analogy, we say that~$\beta$ is \em{split} if $\hat\beta $ is a~split diagram, which means that~$\beta$
contains neither $\sigma_i$ nor $\sigma_i^{-1}$ for some~$i$.

For every link~$L$ there is a~braid $\beta \in B_n$ with $\hat\beta =L$.
The minimal~$n$ for given~$L$ is called \em{braid index} of~$L$.
(See, e.g.,~\cite{WilFr, Morton}.)

\subsection{Burau representation}

The (reduced) \em{Burau representation}~$\psi$ of $B_3$ in $M_2({\mathbb Z}[t^{\pm 1}])$ is def\/ined~by
\begin{gather*}
\psi(\sigma_1) =\left[
\begin{matrix}-t & 1
\\
0 & 1
\end{matrix}
\right],
\qquad
\psi(\sigma_2) =\left[
\begin{matrix}1 & 0
\\
t & -t
\end{matrix}
\right].
\end{gather*}
Then for $k\in{\mathbb Z}$ we have
\begin{gather}
\label{12}
\psi(\sigma_1^k) =\left[
\begin{matrix}(-t)^k & \dfrac{1-(-t)^k}{1+t}
\vspace{1mm}\\
0 & 1
\end{matrix}
\right],
\qquad
\psi(\sigma_2^k) =\left[
\begin{matrix}1 & 0
\vspace{1mm}\\
t \dfrac{1-(-t)^k}{1+t} & (-t)^k
\end{matrix}
\right].
\end{gather}

For a~closed 3-braid, there is a~relation between the Burau representation and the Jones polynomial, which is known from
the related Hecke algebra theory explained in~\cite{Jones}:
\begin{gather}
\label{V_3}
V_{\hat\beta}(t) =\big({-}\sqrt{t}\big)^{[\beta]-2} \left[t\cdot \tr\psi(\beta)+\big(1+t^2\big)\right].
\end{gather}
We will more often than the formula itself use an important consequence of it: for a~3-braid, two of the Burau trace,
exponent sum and Jones polynomial (of the closure) determine the third.

Note that~$\psi$ is faithful on~$B_3$.
One proof which uses directly this relation to the Jones polynomial is given in~\cite{ntriv}.
This property of~$\psi$ is not mandatory for our work, but we use it to save a~bit of exposition overhead at some
places.

More importantly, there is a~way to identify for given matrix whether it is a~Burau matrix, and if so, of which braid.
The Burau matrix determines (for 3-braids) along with the Jones polynomial also the skein and Alexander polynomial.
These in turn determine (the skein polynomial precisely~\cite{3br}, the Alexander polynomial up to a~twofold
ambiguity~\cite{mwf}) the minimal length of a~band representation of the braid.
Thus one has only a~f\/inite list of band representations to check for given Burau matrix.
We used this method in fact also here to identify certain braids from their matrix.
No tools even distantly comfortable are available (or likely even possible) for higher braid groups.

\subsection{Some properties of everywhere equivalent diagrams}\label{S@}

We stipulate that in general~$D$ will be used for a~link diagram and $D'$ for a~diagram obtained from~$D$ by exactly one
crossing change.
If we want to indicate that we switch a~crossing numbered as~$i$, we also write $D'_i$.
Similarly~$\beta$ will stand for a~braid, usually a~particular word of it, and $\beta '$ for a~word obtained by inverting
exactly one letter in~$\beta$.

The central attention of this study is the following type of diagrams.
\begin{defn}
\label{Md}
We call a~link diagram~$D$ \em{everywhere equivalent $($EE$)$} if all diagrams $D'_i$ depict the same link for all~$i$.
(This link may be dif\/ferent from the one represented by~$D$.) If all $D'$ are unknot diagrams, we call~$D$
\em{everywhere trivial}.
\end{defn}

A, potentially complete, list of \em{knotted} everywhere trivial diagrams was determined in~\cite{nikos}.
These are given below and consist of two trefoil and four f\/igure-8-knot diagrams:
\begin{gather}
\label{eeq}
\begin{split}&
\includegraphics{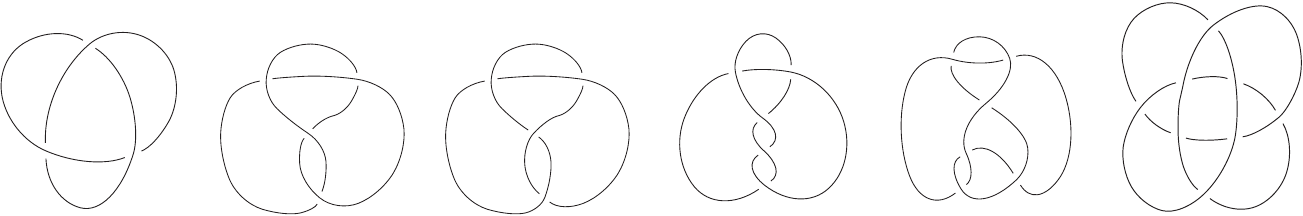}
\end{split}
\end{gather}

We saw this list compatible with the 3-braid case in Theorem~\ref{th1}.
The list draws evidence for its exhaustiveness from various sources.
For minimal crossing diagrams (and in particular that only the trefoil and f\/igure-8-knot occur), the problem seems to
have been noticed previously, and is now believed to be some kind of ``knot-theory folklore''.
Yet, despite its very simple formulation, it is extremely hard to resolve.
Apart from our own (previously quoted) ef\/forts, we are not aware of any recent progress on it.

Another situation where everywhere equivalence can be resolved is for 2-bridge (rational) and Montesinos link diagrams.
The classif\/ication of their underlying links gives a~rather straightforward, albeit somewhat tedious, method to list up
such EE diagrams.
In particular, these diagrams \em{without trivial clasps} are as follows, agreeing also with~\eqref{eeq}: \hfill\null
\begin{wrapfigure}[6]{r}{0.33\textwidth} \vskip8mm
\begin{gather}
\label{xxs}
\begin{split}&
\text{\raisebox{-7.2ex}[0pt][0pt]{\includegraphics{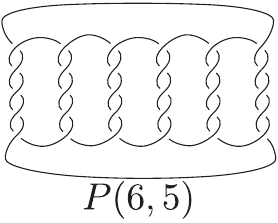}}}
\end{split}
\end{gather}
\end{wrapfigure}
\begin{itemize}\itemsep=0pt\vspace{-8mm}
\item the rational diagrams in~\eqref{eeq} (all except the 8-crossing one),
\item those with Conway notation $C(p)$ and $C(p,-p)$ for $p\ge 1$, and
\item the pretzel diagrams $(\kern1mm\underbrace{\kern-1mm p,\dots,p\kern-1mm}_ {\scbox{$q$ times}}\kern1mm)=P(q,p)$ for
$q\ge 3$ and $p\ge 2$ (see~\eqref{xxs} for an example).
\end{itemize}

Note that in both the Montesinos and 3-braid case, it is not necessary to exclude unknotted everywhere trivial diagrams
to formulate a~reasonable statement.
However, we know that such diagrams occur in multitude outside these two classes, which is another hint to why the
classes are rather special.

The next lemma proposes a~new family of EE diagrams for links, suggested by the 3-braid case.
It identif\/ies the second family in Theorem~\ref{th1} in more general form.

\begin{lem}
When $\beta \in B_n$ is central, then every positive word of~$\beta$ is everywhere equivalent.
\end{lem}

\begin{proof}
The crossing switched versions of $\hat\beta $ are represented up to cyclic permutations by braids of the form
$\sigma_i^{-2}\cdot \alpha$, where~$\alpha$ are cyclically permuted words, one of which represents~$\beta$.
But if~$\beta$ is central, then all~$\alpha$ represent this same central element, and thus all $\sigma_i^{-2}\cdot
\alpha$ are conjugate.
\end{proof}

Beyond this point, we will from now on focus on 3-braids.
The proof of Theorem~\ref{th1} relies on several of their very special properties.
Despite the various input, the length of the argument shows that there was some ef\/fort in putting the pieces together.
In that sense, our initial optimism in carrying out, for example, a~similar investigation of 4-braids, seems little
reasonable.

In relation to our method of proof, we conclude the preliminaries with the following remark.

The use of some algebraic technology might be occasionally (but not always) obsolete, for Birman--Menasco's
work~\cite{BirMen} has reduced the isotopy problem of closed 3-braids (mainly) to conjugacy.
This fact provided a~guideline for our proof.
One place where the analogy surfaces is Lemma~\ref{lm11}, which in some vague sense imitates, on the level of the Jones
polynomial, a~partial case of the combination of Birman--Menasco with the summit power in Garside's conjugacy normal
form~\cite{Garside}.
We will invoke Garside's algorithm at one point quite explicitly, in Lemma~\ref{lm13}.
On the other hand, the use of~\cite{BirMen} would not make the proof much simpler, yet would build a~heavy framework
around it, which we sought to avoid.
We will see that we can manageably work with the Burau representation and Jones polynomial (and that they remain
essential to our argument).

\section{Proof of the non-positive case}\label{S3}

We start now with the proof of Theorem~\ref{th1}.

\subsection{Initial restrictions}

There is also here a~dichotomy as to whether we are allowed to switch crossings of either sign, i.e., whether the
diagram is (up to mirroring) positive or not.
Let us throughout the following call a~braid word everywhere equivalent (EE) if the diagram $\hat\beta $ is such.

We start by dealing with non-positive braids.
The goal is to obtain the f\/irst family in Theorem~\ref{th1}.
Notice here that for any (non-trivial) such word~$\beta$, the closure $\hat\beta $ is either the unknot or the f\/igure-8-knot.

For non-positive braids~$\beta$, a~strong restriction enters immediately, which will play a~central role in the
subsequent calculations.
Since for a~non-positive diagram $D=\hat\beta $, one can get (3-string) braids
of exponent sum dif\/fering by $\pm 4$ representing the same link, it follows from the Morton--Williams--Franks
inequality~\cite{Morton,WilFr} that the skein
(HOMFLY-PT) polynomial~$P$ of such a~link must have a~single non-Alexander variable degree.
Then a~well-known identity~\cite[Proposition~21]{LickMil2} implies that $P=1$.
By~\cite{3br} we can conclude then that \em{$\hat\beta '=D'$ is the unknot}.
In particular, the closure of~$\beta$ is a~knot, and its exponent sum must be zero:
\begin{gather}
\label{btt}
[\beta]=0.
\end{gather}

We remark that if $[\beta ']=\pm 2$ (and the closure is unknotted), then by~\eqref{V_3} its Burau trace is
\begin{gather}
\label{tt}
\tr\psi(\beta ') =(-t)^{\pm 1}.
\end{gather}
Note that not only does this trace determine a~trivial Jones polynomial, but also that the Jones polynomial detects the
unknot for 3-braids (see~\cite{ntriv}).
Thus this trace condition is in fact equivalent to the braid having unknotted closure.

Since we know \em{a priori}
that we expect a~f\/inite answer, it turns out helpful f\/irst to rule out certain subwords of~$\beta$.

Let us f\/irst exclude the cases when~$\beta$ contains up to cyclic permutations subwords of the form $\sigma_i^{\pm
1}\sigma_i^{\mp 1}$, i.e., that it is not cyclically reduced.

It is clear that if such~$\beta$ is everywhere equivalent, then so is the word obtained under the deletion of
$\sigma_i^{\pm 1}\sigma_i^{\mp 1}$.
When we have proved the exhaustiveness of family~\ref{type1} in Theorem~\ref{th1}, we see that it is enough to show that
all words obtained by inserting $\sigma_i^{\pm 1}\sigma_i^{\mp 1}$ cyclically somewhere into any of these words~$\beta$
gives no everywhere equivalent word.

Most cases can be ruled out right away.
Note that in such a~situation for some word~$\beta$ both~$\beta \sigma_i^2$ and~$\beta \sigma_i^{-2}$ must have unknotted closure.
In particular, $[\beta]=0$, and the positive words~$\beta$ need not be treated.

The other (non-positive) words~$\beta$ can be ruled out by noticing that when $\beta \sigma_i^2$ gives (under closure) the
unknot and~$\beta$ the unknot or f\/igure-8-knot, then by the skein relation~\eqref{1} of the Jones polynomial, $\beta
\sigma_i^{-2}$ will have a~closure with some non-trivial polynomial, and so it will not be unknotted.

Thus from now on we assume that~$\beta$ is a~reduced word and has an exponent vector.
The case of exponent vector of length (i.e., weight) $2$ is rather easy, and leads only to $\sigma_1\sigma_2^{-1}$, so
let us exclude this in the following.

\subsection{Syllable types}

We regard thus now~$\beta$ as cyclically reduced, and start by examining what type of syllables can occur in the exponent vector.
For a~syllable, we are interested in the exponent (i.e., sign and length) of the syllable, and the sign of the preceding
and following syllable.

Up to mirroring and $\sigma_1\leftrightarrow\sigma_2$ we may restrict ourselves to positive syllables of $\sigma_2$, and
up to inversion we have three types of signs of neighboring syllables (of $\sigma_1$).

\begin{caselist}
\case Two positive neighboring syllables.
Up to cyclic permutations we have for $n\ge 1$,
\begin{gather}
\label{23}
\beta =\alpha\sigma_1\sigma_2^n\sigma_1.
\end{gather}
Let us call the subword starting after~$\alpha$ the \em{visible subword} of~$\beta$, and its three syllables (at least
two of which are trivial) the \em{visible syllables}.
We cannot assume that~$\alpha$ ends on $\sigma_2^{\pm 1}$, so that the f\/irst visible syllable of~$\beta$ may not be
a~genuine syllable.

We try to f\/ind out now how
$ 
M=\psi(\alpha)$
should look like.
First note that because of~\eqref{btt}, we must have $[\alpha]=-2-n$, and so
\begin{gather}
\label{24}
\det M=(-t)^{-2-n}.
\end{gather}
Next, the condition that switching a~crossing in any of the visible syllables of~$\beta$ should give the unknot implies
from~\eqref{tt} that
\begin{gather*}
\tr\left(M\cdot\psi([-12^n1])\right) = (-t)^{-1},
\qquad
\tr\left(M\cdot\psi([12^{n-2}1])\right) = (-t)^{-1},
\\
\tr\left(M\cdot\psi([12^n-1])\right) = (-t)^{-1}.
\end{gather*}
(Here $(-t)^{-1}$ has to occur everywhere on the right, since we always switch a~positive crossing in~$\beta$
with~\eqref{btt}, and thus $[\beta ']=-2$.)

These three equalities give af\/f\/ine conditions on~$M$, which restrict it onto a~line in the space of $2\times 2$
matrices, regarded over the fraction f\/ield ${\cal F}$ of ${\mathbb Z}[t,t^{-1}]$.
Then the quadratic determinant condition~\eqref{24} will give two solutions.
These live \em{a priori} only in a~quadratic extension of ${\cal F}$.
For the existence of~$\alpha$ we need, (1) that the entries of~$M$ lie outside a~quadratic f\/ield over ${\cal F}$
(i.e.,~the discriminant is a~square), that then (2) they lie in ${\mathbb Z}[t,t^{-1}]$ (i.e., that denominators disappear
up to powers of~$t$), and that in the end (3) the entries build up a~valid Burau matrix of a~braid.
Something that startled us is that we encountered redundant cases which survived until any of these three intermediate
stages.

Putting the equations into MATHEMATICA$^{\rm TM}$
\cite{Wolfram}, using the formulas~\eqref{12}, gives the solutions for
the lower right entry~$d$ of~$M$.
We will use from now on the extra variable
\begin{gather*}
u=-t
\end{gather*}
to simplify the complicated expressions somewhat (even just a~sign af\/fects the presentation considerably!).

The value for~$d$ determined by MATHEMATICA looks thus:
\begin{gather}
d=\bigg[u^{7+n}+5u^{3+2n}-7u^{4+2n}+6u^{5+2n}-5u^{6+2n}+4u^{7+2n}-2u^{8+2n}-3u^{2+3n}
\nonumber
\\
\phantom{d=}
{}+8u^{3+3n}-10u^{4+3n}+6u^{5+3n}-2u^{6+3n}+u^{7+3n}+u^{3+4n}+u^{5+n}t+u^{1+2n}t+u^{3+4n}t
\nonumber
\\
\phantom{d=}
{}\pm \Big(u^{1+2n}{(1+t)}^2 \big(4u^n+8u^{2+n}+8u^{4+n}-4u^{5+n}+4u^{6+n}+4u^nt
\label{horror}
\\
\phantom{d=}
{}+ \big(3+10u^n+3u^{2n}\big) t^3 \big)
 {\big(u^n+3u^{2+n} +(2+u^n) u^nt+u^nt^3+(-1+u^n)t^4 \big)}^2\Big)^{1/2}\bigg]
\nonumber
\\
\phantom{d=}
{}\times \Bigl({2u^{2n}t(1+t)^2\big[u^n+3u^{2+n}+3u^{4+n}+u^{6+n}+u^nt+ \bigl(1+4u^n+u^{2n} \bigr) t^3+u^nt^5 \big]}\Bigr)^{-1}.
\nonumber
\end{gather}
Our attention is dedicated f\/irst to the discriminant (occurring under the root).
Removing obvious quadratic factors, we are left with the polynomial
\begin{gather}
\label{??}
u\big({-}3u^3+4u^n -4u^{n+1}+ 8u^{2+n}-10u^{n+3}+8u^{4+n}-4u^{5+n}+ 4u^{6+n} -3u^{2n+3}\big).
\end{gather}
We need that this polynomial becomes the square of a~polynomial in ${\mathbb Z}[t,t^{-1}]$.

For $n\ge 4$ the edge coef\/f\/icients are $-3$, and thus the polynomial is not a~square.
Similarly for $n=2$, where the minimal and maximal degrees are odd.

For $n=1$ the visible subwords we consider in~\eqref{23} are just the half-center element~$\Delta$ (and, under various
symmetries, its inverse).
Their exclusion as subwords of~$\beta$ will be most important for the rest of the argument (see Lemma~\ref{lm13}), but
occurs here in the most peculiar way.

We have in~\eqref{??} the polynomial $4-4u+5u^2-10u^3+5u^4-4u^5+4u^6$.
This becomes a~square, and so we obtain values for~$d${\samepage
\begin{gather*}
d=\frac{t^4+2t^5-t^6-4t^7-t^8+2t^9+t^{10}\pm{\sqrt{t^6{(1+t)}^8\left(-2+t-4t^2+t^3-2t^4 \right)^2}}}{2t^4 {(1+t)}^4 \left(1-t+3t^2-t^3+t^4 \right)}
\end{gather*}
in ${\cal F}$ (rather than some of its quadratic f\/ields).}

For the choice of negative sign we get for~$d$
\begin{gather*}
d =\frac{1+t^2+t^4}{t-t^2+3t^3-t^4+t^5}.
\end{gather*}
The evidence that this is not a~Laurent polynomial in~$t$ can be sealed, for example, by setting $t=\frac{1}{2}$.
The expression evaluates to $\frac{42}{19}$, whose denominator is not a~power of~2.

For the `+' we get $d=-1/t$.
This gives indeed a~matrix in ${\mathbb Z}[t,t^{-1}]$:
\begin{gather*}
M =\left[
\begin{matrix}
t^{-2} & 0
\\
t^{-1} & -t^{-1}
\end{matrix}
\right].
\end{gather*}
But there is no braid~$\alpha$ with such Burau matrix.
This can be checked via the Alexander polynomial of the (prospective) closure, but there is a~direct argument (which
appeals, though, to the faithfulness of~$\psi$).
The matrix in question is $M=t^{-2}\cdot \psi(\sigma_2)$, but scalar Burau matrices are only in the image of the center
of $B_3$, and these are powers of $t^3$.

For $n=3$ the polynomial~\eqref{??} is $1-4u+8u^2-10u^3+8u^4-4u^5+u^6$, which is a~square, and the rather complicated
expressions~\eqref{horror} become
\begin{gather*}
\frac{t^8+4t^9+6t^{10}+3t^{11}-3t^{12}-6t^{13}-4t^{14}-t^{15}\pm{\sqrt{t^{16}{(1+t)}^{10}{(1+t+t^2)}^2}}}{2t^{11} {(1+t)}^4(1+t+t^2)},
\end{gather*}
giving $d=t^{-3}$ and $d=-t^{-2}$.
These lead to the matrices
\begin{gather*}
M =\left[
\begin{matrix} -t^{-2} & t^{-3}
\\
0 & t^{-3}
\end{matrix}
\right]
\qquad
\text{and}
\qquad
M =\left[
\begin{matrix} t^{-3} & 0
\\
t^{-2} & -t^{-2}
\end{matrix}
\right].
\end{gather*}
We used the Alexander polynomial to check that these indeed occur as Burau matrices, and to see what their braids are.
(We remind that~$\psi$ is faithful on $B_3$.) The answer is
\begin{gather*}
\alpha_1=\Delta^{-2}\sigma_1
\qquad
\text{and}
\qquad
\alpha_2=\Delta^{-2}\sigma_2.
\end{gather*}
These solutions were unexpected, but can be easily justif\/ied: putting $\alpha_i$ for~$\alpha$ in~\eqref{23} for $n=3$,
one easily sees that switching any of the last 5 crossings gives a~braid with unknotted closure.

\case One positive and one negative neighboring syllable.
In this case we have
\begin{gather}
\label{23'}
\beta =\alpha\sigma_1^{-1}\sigma_2^n\sigma_1.
\end{gather}
Thus $[\alpha]=-n$, and we have to check the identity for~$M$ as in~\eqref{24}
\begin{gather*}
\det M=(-t)^{-n},
\end{gather*}
and the unknot trace conditions:
\begin{gather*}
\tr\left(M\cdot\psi([12^n1])\right) = -t,
\qquad
\tr\left(M\cdot\psi([-12^{n-2}1])\right) = (-t)^{-1},
\\
\tr\left(M\cdot\psi([-12^n-1])\right) = (-t)^{-1}.
\end{gather*}
(Now in the f\/irst case we switch a~negative crossing, thus $[\beta ']=2$, and we need the trace $-t$.)

The solutions for the lower right entry~$d$ of~$M$ are now even less pleasant than~\eqref{horror}, and thus we do not
reproduce them here.
It is again important mainly to look at the discriminant, which was identif\/ied by MATHEMATICA as
\begin{gather*}
\begin{split}&
 u^{1+2n}{(1+t)}^2 {\left(1-u^n+t+t^2 \right)}^2 \Bigl(-4u^{2n}+4 \left(1+u^n
+ u^{2n} \right) u^{5+n}-12u^{2+2n}-12u^{4+2n}
\\
&\qquad
{} -4u^{6+2n}-4\left(1+u^n+u^{2n} \right) u^nt+\left(3-4u^n-6u^{2n}-4u^{3n}+3u^{4n} \right) t^3 \Bigr),
\end{split}
\end{gather*}
and becomes a~square times the following polynomial:
\begin{gather*}
u \bigl({-}3u^{4n+3} +4u^{3n+5} +4u^{3n+3} +4u^{3n+1} -4u^{2n+6} +4u^{2n+5} -12u^{2n+4}+\\
\qquad
{}+6u^{2n+3} -12u^{2n+2} +4u^{2n+1} -4u^{2n} +4u^{n+5} +4u^{n+3} +4u^{n+1} -3u^3\bigr).
\end{gather*}
This polynomial is not a~square for $n\ge 3$ by the leading coef\/f\/icient argument.
It is, though, a~square for $n=1,2$, and gives solutions for~$d$,~$M$, and~$\alpha$.
The braids~$\alpha$ are more easily identif\/ied by direct observation.
For $n=2$ we have
\begin{gather*}
\alpha_1=\sigma_2^{-1}\sigma_1^{-1}
\qquad
\text{and}
\qquad
\alpha_2=\sigma_1^{-1}\sigma_2^{-1}.
\end{gather*}
The solution $\alpha_1$ can be seen from the word $[-2-1-1221]$ in family~\ref{type1}.
The other solution (was guessed but) is also easily conf\/irmed.
For $n=1$ both solutions stem from words in family~\ref{type1}:
\begin{gather*}
\alpha_1=\sigma_2^{-1}
\qquad
\text{and}
\qquad
\alpha_2=\sigma_2^{-1}\sigma_1^{-1}\sigma_2\sigma_1\sigma_2^{-1}.
\end{gather*}

\case Two negative neighboring syllables.
Here we have with
\begin{gather}
\label{23''}
\beta =\alpha\sigma_1^{-1}\sigma_2^n\sigma_1^{-1},
\end{gather}
the condition $\det M=(-t)^{2-n}$, and the trace conditions
\begin{gather*}
\tr\left(M\cdot\psi([12^n-1])\right) = -t,
\qquad
\tr\left(M\cdot\psi([-12^{n-2}-1])\right) = (-t)^{-1},
\\
\tr\left(M\cdot\psi([-12^n1])\right) = -t.
\end{gather*}
The solutions for the lower right entry of~$M$ look thus:
\begin{gather*}
d= -\Bigl(2u^{2+n}-2u^{3+n}+u^{4+n}+3u^{2+2n}-4u^{3+2n}+u^{4+2n}+u^{2+3n}
+ u^nt+u^{3n}t
\\
\phantom{d=}
 {}\pm {\sqrt{u^{1+2n}{(1+t)}^2 \left(4u^n+4u^{2+n}+\left(3+2u^n+3u^{2n}
\right) t \right)  {\left(1-u^n+t+t^2 \right)}^2}} \Bigr)
\\
\phantom{d=}
 {} \times \frac{1}{2u^{2n} {(1+t)}^2\left(1+t+u^nt+t^2 \right)}.
\end{gather*}
The non-square part of the discriminant is
\begin{gather*}
-3u^2+4u^{1+n}-2u^{2+n}+4u^{3+n}-3u^{2+2n}.
\end{gather*}
Here the edge coef\/f\/icients become $-3$ for $n\ge 2$.
Thus $n=1$.
Now~$M$ again becomes a~Burau matrix, and there are the solutions
\begin{gather*}
\alpha_1=\sigma_1
\qquad
\text{and}
\qquad
\alpha_2=\sigma_2.
\end{gather*}
(Again the f\/irst comes from $[1-21-2]$ in family~\ref{type1}.)
\end{caselist}

With the discussion in the preceding three cases, we have thus now obtained restrictions on how syllables in~$\beta$ can
look like.
There are four `local' syllable types (up to symmetries), which can be summarized thus, and we will call
\em{admissible}.
\begin{itemize}\itemsep=0pt
\item No syllable has length at least 4.
\item A~syllable of length 3 has both neighbored syllables of the same sign.
\item A~syllable of length 2 has exactly one of its two neighbored syllables having the same sign.
\item A~syllable of length 1 has at most one of its two neighbored syllables having the same sign.
\end{itemize}
For each admissible syllable, we have also identif\/ied the two possible braids outside the syllable (although not their
word presentations).

\subsection{Words of small length}

The next stage of the work consists in verifying a~number of words of small length.
One can easily argue (see~\cite{e12cp}) that for more than one component, there are no crossings in an everywhere
equivalent diagram between a~component and itself.
We notice that this observation specializes here to saying that (for non-split braids) the exponent sum (or word length)
is even.
This will not be essential in the proof, but helpful to avoid checking certain low crossing cases.

The test of small length words is done by an algorithm which goes as follows.

We start building a~word~$\beta$ of the form
\begin{gather*}
\beta =\alpha\gamma,
\end{gather*}
where $\gamma$ is a~word known to us, and we know the Burau matrix
\begin{gather*}
M=\psi(\alpha)
\end{gather*}
of~$\alpha$.
(Thus essentially we know~$\alpha$ itself, we do not know a~word presentation of it.) We perform this with each of the
shapes in~\eqref{23},~\eqref{23'} and~\eqref{23''} with the two possible values of~$\alpha$ we found in each case.
(Thus there are in total 8 initial pairs of parameters $(M,\gamma)$.) The understanding is that switching any of the
crossings in $\gamma$ gives a~braid with unknotted closure.
We call $\gamma$ an \em{extendable word}, in the sense that it can potentially be extended to a~solution~$\beta$.

Whenever~$M$ is the identity matrix, we can take $\beta =\gamma$ and have an everywhere equivalent braid, which we output.

Next we try to extend $\gamma$ by one letter $\tau=\sigma_i^{\pm 1}$, so that it is not the inverse of the preceding
letter, and the admissible syllable shapes are not violated.
Let ${\tilde M}=\psi(\tau^{-1}) \cdot M$.
We test whether
\begin{gather*}
\tr\big({\tilde M}\cdot\psi\big(\gamma\cdot \tau^{-1}\big)\big)=(-t)^{\mp 1},
\end{gather*}
which is equivalent to whether a~crossing change at the new crossing (also) gives the unknot.
If this happens, we can continue the algorithm with $\gamma$ replaced by $\gamma\cdot\tau$ and~$M$ replaced by ${\tilde M}$.

This procedure can yield the solutions~$\beta$ up to given number of crossings (word length), and also produce the list
of extendable braids $\gamma$ up to that crossing number.

Note that, since potential EE solutions can be directly checked, it is often enough to work with particular values of~$t$.
(These can be \em{a priori} complex numbers, but in practice most helpfully should be chosen to be rational.) We did not
feel conf\/ident about this in the three initial cases, because of the presence of variable exponents.
Alternatively, we could have used a~$t$ with $|t|<1$ and some convergence (and error estimation) argument, but whether
that would have made the proof nicer is doubtful.

For rational~$t$, we implemented the above outlined procedure in C++, whose arithmetic is indef\/initely faster than
MATHEMATICA.
It, however, has the disadvantage of thoughtlessly producing over-/underf\/lows, and some ef\/fort was needed to take care
of that and to work with rational numbers whose numerator and denominator are exceedingly large.

With this problem in mind, it is recommended to use simple (but non-trivial) values of~$t$.
We often chose $t=2$, but also $t=3$, and a~few more exotic ones like $t=4/5$ whenever feasible.
We were able to perform the test up to 15 crossings for $t=2$ and up to 12 crossings (still well enough for what we
need, as we will see below) for the other~$t$.

This yielded the desired family~\ref{type1}, but also still a~long list of extendable words even for large crossing number.
Such words were not entirely unexpected: one can see, for example, that when~$\beta$ is a~solution, then any subword
$\gamma$ of a~power $\beta^k$ of~$\beta$ is extendable (with~$\alpha$ being, roughly, a~subword of $\beta^{1-k}$).
There turned out to be, however, many more extendable words, which made extra treatment necessary.

The following argument gives a~mild further restriction.

Assume~$\beta$ has a~subword $[12-1-2]$.
Switching the `$2$' will give $[-2-1]$, while switching the~`$-1$' will give $[21]$.
Now, these two subwords can be realized also by switching either crossings in $[-21]$, which is a~word of the same braid as $[12-1-2]$.
This means that if~$\beta$ is EE, then also a~word is EE in which $[12-1-2]$ was replaced by $[-21]$.
Thus verifying words up to 10 crossings would allow us to inductively discard words containing $[12-1-2]$ and its various equivalents.

\subsection{Global conditions}

All these `local' conditions were still not enough to rule out all possibilities, and f\/inally we had to invent a~`global' argument.

For this we use the following fact: if $\beta '\in B_3$ has unknotted closure, then $\beta '$ is conjugate to $\sigma_1^{\pm 1}\sigma_2^{\pm 1}$.
The f\/irst proof appears to be due to Murasugi~\cite{Murasugi}.
It was recovered by Birman--Menasco~\cite{BirMen}.
A~dif\/ferent proof, based on the Jones polynomial, follows from (though not explicitly stated in)~\cite{ntriv}.
Namely, one can use relations of the sort
\begin{gather}
\label{o9}
\sigma_1\sigma_2^k\sigma_1^{-1} =\sigma_2^{-1}\sigma_1^k\sigma_2,
\end{gather}
which do not augment word length (together with their versions under the various involutions), and cyclic permutations
to reduce $\beta '$ to a~length-2 word.

The non-conjugacy to $\sigma_1^{\pm 1}\sigma_2^{\pm 1}$ of a~crossing-switched version $\beta '$ of our braids~$\beta$ is
detected by Garside's algorithm~\cite{Garside}.
We adapt it in our situation as follows.

\begin{lem}
\label{1111}
Assume a~braid $\beta '\in B_3$ is written as $\Delta^k\alpha$ with~$\alpha$ a~positive word with cyclically no trivial
syllables and $k\le -2$.
Then $\beta '$ is not conjugate to $\sigma_1^{\pm 1}\sigma_2^{\pm 1}$.
\end{lem}

\begin{proof}
This is a~consequence of Garside's summit power in the conjugacy normal form.
There is again an alternative (but a~bit longer) way using $\spaan_t\tr\psi(\beta ')$.
We refer for comparison to the proof of Lemma~\ref{lm11}, but for space reasons only brief\/ly sketch the argument here.
For~$\alpha$ one uses the relation~\eqref{V_3} and that the span of $V(\hat\alpha)$ is determined by the adequacy of the diagram~$\hat\alpha$.
The center multiplies~$\tr\psi$ only by powers of~$t^3$.
\end{proof}

\begin{lem}
\label{lm13}
Let $\beta '$ be a~reduced braid word $($not up to cyclic permutations$)$ with no~$\Delta$ and~$\Delta^{-1}$ as subwords.
Then under braid relations $($without cyclic permutations$)$ one can write
\begin{gather*}
\beta '=\Delta^{-k}\alpha
\end{gather*}
with $k\ge 0$ and~$\alpha$ a~positive word.
Moreover,
\begin{itemize}
\itemsep=0pt
\item~$k$ is at least half of the number of negative letters in $\beta '$.
\item Only the first and last syllable of~$\alpha$ can be trivial.
If $\beta '$ starts and ends with positive letters~$\tau$ $($not necessarily the same for start and end$)$, which are not
followed resp.\ preceded by $\bar\tau$, then~$\alpha$ has no trivial syllable.
\end{itemize}
\end{lem}
(For the bar notation recall the end of Section~\ref{22}.)

\begin{proof}
This is the result of the application of Garside's procedure on $\beta '$.
We manage, starting with trivial $\alpha'$ and $k'=0$, a~word presentation of $\beta '$
\begin{gather*}
\Delta^{-k'}\alpha'\gamma,
\end{gather*}
with the following property: $\alpha'$ is a~positive word with only the f\/irst and last syllable possibly trivial, and
$\gamma$ is a~terminal subword of $\beta'$.
We also demand that if $\gamma$ starts with a~positive letter, this is the same as the f\/inal letter of $\alpha'$ (unless $\alpha'$ is trivial).
We call this the \em{edge condition}.

We apply the following iteration.
\begin{enumerate}
\itemsep=0pt
\item
\label{step1}
Move as many positive initial letters from $\gamma$ as possible into the end of $\alpha'$, so as $\gamma$ to start with
a~negative letter~$\tau$.
This will not produce internal trivial syllables in $\alpha'$ because of the edge condition and because $\beta '$ has
no~$\Delta$ subword.

\item If $\gamma$ has no negative letters, then we move it out into $\alpha'$ entirely, and are done with $k=k'$.

\item Otherwise we apply two types of modif\/ications:
\begin{itemize}
\itemsep=0pt
\item If $\gamma$'s second letter (exists and) is $\bar\tau$, delete $\tau\bar\tau$ in $\gamma$, replace $\alpha'$~by
$\bar\alpha'$, add at its end~$\tau^{-1}$, and augment~$k$ by~$1$.
\item If $\gamma$'s second letter (does not exist or) is not $\bar\tau$, delete~$\tau$ in $\gamma$, replace $\alpha'$~by
$\bar\alpha'$, add at its end $\tau^{-1}\bar\tau^{-1}$, and augment~$k$ by~$1$.
\end{itemize}

Then go back to step~\ref{step1}.
\end{enumerate}

In the end we obtain the desired form.
Since there is no $\Delta^{-1}$ in $\beta '$, each copy of $\Delta^{-1}$ added compensates for at most two negative
letters of $\beta '$.
\end{proof}

Now we consider an EE word~$\beta$ of, say, more than 10 crossings.
We switch, in a~way we specify below, a~crossing properly and apply the above procedure starting at a~cyclically
well-chosen point of the resulting word $\beta '$.
We obtain the shape of Lemma~\ref{lm13}.
This gives the contradiction that the closure $\hat\beta '$ is not unknotted by Lemma~\ref{1111}.

If~$\beta$ contains a~syllable of length 2 or 3, then we have (up to symmetries) $[12221]$ or $[-1221]$.
Switch the f\/inal~`$1$'.
This does not create a~$\Delta^{-1}$ because there is no subword $[21-2-1]$ in~$\beta$.
Then apply the procedure starting from the second last letter: $\beta '=[2-1\cdots-12]$ or $[2-1\cdots 122]$.

With this we can exclude non-trivial syllables in~$\beta$.
Next, if there is a~trivial syllable between such of opposite sign, up to symmetries $[-2-12]$, we have with its two
further neighbored letters
\begin{gather}
\label{q1}
[1|-2-12|-1].
\end{gather}
The right neighbor cannot be `$1$', because we excluded $[-2-121]$ as subword.
The left neighbor cannot be `$-1$', because we excluded $\Delta^{-1}$.
Now we switch the middle `$-1$' in the portion~\eqref{q1}, and start the procedure with the following (second last in
that presentation) letter `$2$'.

The only words that remain now are the alternating words $\beta = [(1-2)^k]$.
There are various ways to see that they no longer create problems.
In our context, one can switch a~positive letter, group out a~$\Delta^{-1}$ built together with the neighbored letters,
and then start the procedure right after that $\Delta^{-1}$.

This completes the proof of the non-positive braids.

\section{Proof of the positive case}\label{S4}

\subsection{Adequate words}

{}From now we start examining, and gradually excluding the undesired, positive braids in $B_3$.
The nature of this part is somewhat dif\/ferent.
Here no electronic computations are necessary, but instead a~delicate induction argument.

The presence of the central braids and their realization of every positive word as subword explain that no `local'
argument can work as in the non-positive case.
Thus from the beginning we must use the `global' features of the braid words.

Our attitude will be that except for the stated words~$\beta$, we f\/ind two diagrams $D'=\hat\beta '$ we can distinguish~by
the Jones polynomial.
Because of the skein relation~\eqref{1} of~$V$, one can either distinguish the Jones polynomial of (the closure of) two
properly chosen crossing-switched versions $\beta '$, or of two smoothings of~$\beta$.
(In the crossing-switched versions, a~letter of~$\beta$ is turned into its inverse, and in the smoothings it is deleted.)

Moreover, one can switch back and forth between the Jones polynomial and the Burau trace, because of the consequence
of~\eqref{V_3} stated below it.
Notice that for a~positive word, length and exponent sum are the same.

Accordingly we call a~word \em{$\psi$-everywhere equivalent}, if all crossing-switched versions (or equivalently, all
smoothed versions) have the same Burau trace.

Trivial syllables will require a~lot of attention in the following, and thus, to simplify language, we set up the
following terminology.

\begin{defn}
We call a~positive word \em{adequate}, resp.\ \em{cyclically adequate}, if it has no trivial syllable (the exponent
vector has no `1'), resp.\ has no such syllable after cyclic permutations.
Otherwise, the word is called \em{$($cyclically$)$ non-adequate}.
\end{defn}

This choice of language is suggested by observing that cyclically adequate words~$\beta$ give adequate diagrams $\hat\beta$.
Contrarily, for cyclically non-adequate words~$\beta$ the diagrams $\hat\beta $ are not adequate: a~trivial syllable of
a~positive word~$\beta$ always gives a~self-trace in the~$B$-state of $\hat\beta $, i.e., $\hat\beta $ is not~$B$-adequate.
(However, being positive, $\hat\beta $ is always~$A$-adequate, and thus it is \em{not} inadequate in the sense
of~\eqref{aqdf}.)

A useful application of adequacy is the following key lemma, which will help us carry out the induction without
unmanageable calculations.
We call below a~trivial syllable \em{isolated} if it is not cyclically followed or preceded by another trivial syllable.
Recall also the weight $\omega(\beta)$ from Section~\ref{22}.

\begin{lem}
\label{lm11}
Assume positive words~$\beta$ and $\gamma$ have the same length, and $\gamma$ is cyclically adequate $($i.e., its cyclic exponent vector has no {\rm `1')}.
Further assume that either
\def\labelenumi{{\rm \theenumi)}}
\begin{enumerate}\itemsep=0pt
\item $\omega(\gamma)<\omega(\beta)$, or
\item $\omega(\gamma)=\omega(\beta)$ and~$\beta$ has an isolated $($trivial$)$ syllable.
\end{enumerate}
Then $\tr\psi(\beta)\ne \tr\psi(\gamma)$.
\end{lem}

\begin{proof}
It is enough to argue with the Jones polynomial.
The closure diagram $\hat\gamma$ is adequate, and by counting loops in the~$A$- and~$B$-states, we see
\begin{gather}
\label{111}
\spn V(\hat\gamma)=[\gamma]-\omega(\gamma)+1.
\end{gather}
Similarly, by comparing the extremal degrees in the bracket expansion~\eqref{eq_12} of the diagram $\hat\beta $, we have
\begin{gather}
\label{112}
\spn V(\hat\beta)\le [\beta]-\omega(\beta)+1,
\end{gather}
and $[\gamma]=[\beta]$.
The f\/irst case of the lemma is then clear from~\eqref{111} and~\eqref{112}.

If $\omega(\gamma)=\omega(\beta)$, the right hand-sides of~\eqref{111} and~\eqref{112} agree, so we argue that the
inequali\-ty~\eqref{112} is strict.
When a~trivial syllable in~$\beta$ is isolated, so is its (self-)trace in the~$B$-state of $\hat\beta $, as def\/ined below~\eqref{lpt}.
In follows then from the explained work of~\cite{BaeMor} that the extreme~$B$-degree term is zero, making~\eqref{112} strict.
\end{proof}

We use the following lemma to f\/irst get disposed of cyclically adequate braids.
Let us from now on use the symbol `$\doteq$' for \em{equality of braid words up to cyclic permutations}.

\begin{lem}
\label{lm0}
Assume~$\beta$ is a~positive adequate~$\psi$-everywhere equivalent word.
Then all cyclic exponent vector entries are equal, i.e., $\beta \doteq[(1^l2^l)^k]$.
\end{lem}

\begin{proof}
Let~$m$ be the minimal exponent vector entry, and assume by contradiction that there is an exponent vector entry $m'>m$.

If $m=2$, compare the smoothings in the syllables corresponding to~$m$ and $m'$.
There is a~contradiction by Lemma~\ref{lm11}.
For $m=3$, consider the crossing-switched versions.

Now, when $m>3$, one can argue as follows.
Because of the skein relation~\eqref{1}, the Jones polynomial of $\alpha\sigma_i^{m-k}$ for f\/ixed~$k$ is determined~by
the one of $\alpha\sigma_i^{m}$ and $\alpha\sigma_i^{m-1}$:
\begin{gather*}
V\big(\widehat{\alpha\sigma_i^{m-k}}\big) =A_k(t) V\big(\widehat{\alpha\sigma_i^{m}}\big)+B_k(t)V\big(\widehat{\alpha\sigma_i^{m-1}}\big),
\end{gather*}
with $A_k,B_k\in{\mathbb Z}[t^{1/2},t^{-1/2}]$ independent of~$\alpha$ and~$i$.
We apply this argument for $k=m-1$ once on the syllable $\sigma_i^m$ and once on $\sigma_j^{m'}$.
Then we see again two positive words of equal length and weight that must have the same Jones polynomial, one of which
has a~(single, and thus isolated) trivial syllable, and the other has none.
As before, Lemma~\ref{lm11} gives a~contradiction.
\end{proof}

For the rest of the proof, we consider a~cyclically non-adequate word~$\beta$, and use induction over the word length.
We assume that~$\psi$-everywhere equivalent braids of smaller length are in families~\ref{type2}, \ref{type3} and~\ref{type4}.
It will be helpful to make the families disjoint by excluding the (central) cases of $l=1$ and $3\mid k$ in
family~\ref{type2}.

\begin{lem}
\label{lm1}
When~$\beta$ is positive and~$\psi$-everywhere equivalent and~$\beta$ has a~$6$-letter subword representing~$\Delta^2$, then
deleting that subword gives a~$\psi$-everywhere equivalent braid word.
\end{lem}

\begin{proof}
All the crossing switched versions of~$\beta$ outside the copy of $\Delta^2$ have the same Burau trace, and deleting that
copy of $\Delta^2$, the Burau trace multiplies by $t^{-3}$.
\end{proof}

In relation, let us remark that $\Delta^2$ has the following 6-letter presentations up to
$\sigma_1\leftrightarrow\sigma_2$ and inversion (but \em{not} cyclic permutations):
\begin{gather}
\label{Dlw}
[121212],
\qquad
[212212],
\qquad
[211211].
\end{gather}

\begin{lem}
\label{lm2}
If~$\beta$ contains up to cyclic permutations two disjoint subwords~$\Delta$ $($given as $[121]$ or $[212])$,
\begin{gather*}
\beta =\alpha_1\Delta\alpha_2\Delta
\end{gather*}
and~$\beta$ is~$\psi$-everywhere equivalent, then so is $\bar\alpha_1\alpha_2$ $($where the over-bar again denotes
$\sigma_1\leftrightarrow\sigma_2)$.
\end{lem}

\begin{proof}
Note that the crossing changes outside the two copies of~$\Delta$ will commute with putting the two~$\Delta$ close using
\begin{gather}
\label{sld}
\alpha_2\Delta=\Delta\bar\alpha_2
\end{gather}
to form a~$\Delta^2$, and then apply Lemma~\ref{lm1}.
\end{proof}

The move~\eqref{sld} will be used extensively below, and will be called \em{sliding}.
Obviously, one can slide any copy of~$\Delta$ through any subword $\alpha_i$.

\subsection{Induction for words with trivial syllables}

\begin{caselist}

\case
$\beta$ has a~$\Delta^2$ subword.
We can apply Lemma~\ref{lm1}, and use induction.

We have a~central $\Delta^2$ word inserted somewhere in a~\em{remainder}, which is either (a) a~central word, (b)
a~split word $\sigma_1^k$ or (c) a~symmetric word $[1^l2^l]^k$.
In (a) we have a~central word~$\beta$, and this case is clear.
So consider (b) and (c).

\begin{caselist}
\case
\label{c1q}
The split remainders.
Then $\beta =\Delta^2\sigma_1^k$ for $k>0$.
So we have the words\footnote{In the following boldfaced letters in braid words should be understood as such whose
switch (resp.\ smoothing) yields a~word that can take the function of $\gamma$ in Lemma~\ref{lm11}.
Italic letters switched (resp.\ removed) should yield the corresponding~$\beta$ in that lemma.} (up to reversal, cyclic permutations and
$\sigma_1\leftrightarrow\sigma_2$)
\begin{gather*}
\def\b#1{{\mathbf #1}} \big[|1211\b 21|1^k\big],
\qquad
\big[|1\b 21212|1^k\big],
\qquad
\big[|2\b 12212|1^k\big],
\qquad
\big[|2212\b 21|1^k\big],
\qquad
\big[|\b 211211|1^k\big].
\end{gather*}
We distinguish them directly by two smoothings: for the letters indicated in bold, we obtain a~$(2,n)$-torus link, and
smoothing a~letter in $1^k$ gives a~positive braid with $\Delta^2$, and thus under closure a~link of braid index 3.

\case The symmetric remainders $[1^l2^l]^k$.

\begin{caselist}
\case Let f\/irst $l=1$ (and $3\nmid k$).
If $k>3$, then however one inserts that $\Delta^2$, one has outside (in the remainder) up to cyclic permutations the
word $\delta=[121212]$ or $[212121]$, i.e.,
\begin{gather*}
\beta \doteq \Delta^2\delta\alpha.
\end{gather*}
(For `$\doteq$' recall above Lemma~\ref{lm0}.
Of course,~$\delta$ represents $\Delta^2$, but we separate both symbols as subwords of~$\beta$.) If the whole word~$\beta$
is~$\psi$-everywhere equivalent, then by Lemma~\ref{lm1} (since~$\delta$ represents $\Delta^2$) so is $\Delta^2\alpha$,
obtained after deleting~$\delta$ in the remainder.
Note for later, when we insert~$\delta$ back, that this insertion must be done so that the last letter of~$\delta$ and
the f\/irst letter of~$\alpha$ are not the same.

Iterating this argument, we can start by testing $\alpha=[12]$ and $\alpha=[1212]$.
Thus it is enough to see that when one inserts a~$\Delta^2$ word into (or before or after) $[12]$ or $[1212]$, the
result~$\beta$ is not~$\psi$-everywhere equivalent, unless it is $\beta \doteq [12]^{k}$ for $k=4$ or~$5$.
These are all knot diagrams of $\le 10$ crossings, and they can be checked directly.
Then (iteratedly) inserting back~$\delta$ must be done so as to yield only $\beta \doteq [12]^{k}$ for higher~$k$.

\case
\label{c122_}
Consider next $l\ge 2$.
Then~$\beta$ is up to cyclic permutations of the form $\Delta^2\alpha$ with the exponent vector of~$\alpha$ having
a~`$1$' possibly only at the start and/or end.

We want to apply Lemma~\ref{lm11} to exclude these cases.

\begin{caselist}
\case Let~$\alpha$ be adequate (i.e., have no trivial syllable even at its start and end).
We compare two crossing-switched versions of $\beta =\Delta^2\alpha$.
First we switch a~crossing in $\Delta^2$, turning~$\beta$ into $\sigma_1^2\sigma_2^2\alpha$ or
$\sigma_2^2\sigma_1^2\alpha$, which is a~cyclically adequate word.
Another time we switch a~crossing in (any non-trivial syllable of)~$\alpha$, yielding $\Delta^2\alpha'$ for a~positive
$\alpha'$, whereby the weight decreases by~2.\looseness=-1

Let us make this argument more precise.
We look at the three words for $\Delta^2$ in~\eqref{Dlw}.

{\def\b#1{{\mathbf #1}} For $[1\b 21212]$, switching the boldfaced `2' gives $[1^22^2]$, and the weight decreases by 4.
Thus we are done by the f\/irst (unequal weight) case of Lemma~\ref{lm11}.

For $[2\b 122\b 12]$, switching either one of the boldfaced letters gives $[1^22^2]$ or $[2^21^2]$.
At least for one of these two words the weight decreases by 4 (and we are through as above), unless both neighboring
letters of $\Delta^2=[212212]$ are `2': $[2|212212|2]$.
In this case the weight decreases by 2, as it does for $\Delta^2\alpha'$.
However, in $\Delta^2\alpha'$, at least one of the two `1' of $\Delta^2$ gives an isolated trivial syllable.
Thus we can apply the second case of Lemma~\ref{lm11}.

For $[2\b 11211]$, we switch to $[2^21^2]$, and the weight decreases by~$2$, but $\Delta^2$ in $\Delta^2\alpha'$ has an
isolated trivial syllable `2'.}

\case Thus let~$\alpha$ have a~trivial initial or f\/inal syllable.
Then, because of $l\ge 2$, the initial and f\/inal syllable have the same index (1 or 2).

Up to $\sigma_1\leftrightarrow\sigma_2$, cyclic permutation and inversion, it is enough to look at the cases
\begin{gather}
\label{44}
\def\a#1{\protect{\mathbf #1}}\def\b#1{\protect{\mathit #1}} \catcode`\|=12
\begin{array}
{*3{@{}l@{
\qquad
}}l} \catcode`\|=12 [1|\b 12\a 1121|1], & \catcode`\|=12 [1|\b 12\a 1212|1], & \catcode`\|=12 [1|\a 21\b 2121|1], &
\catcode`\|=12 [1|\a 21\b 2212|1],
\\[1mm] \catcode`\|=12 [1|11\a 2\b 112|12], & \catcode`\|=12 [1|\b 1221\a 22|12], & \catcode`\|=12 [1|\a 22122\b 1|12], &
\catcode`\|=12 [1|\a 2112\b 11|12].
\end{array}
\end{gather}

These words are again dealt with by Lemma~\ref{lm11}.
In all eight cases switching the boldfaced letter gives, after braid relations, a~positive cyclically adequate word, and
turning the italic one gives a~positive cyclically non-adequate word.
To see this, it is helpful to recall and use the identities~\eqref{o9}.

We exemplarily show the fourth word.
Switching the italic letter gives $[1|2112|1]$, with weight decreasing by~$2$.
Switching the boldfaced letter gives
\begin{gather*}
[1\underline{-212}2121]\to[112\underline{-121}21]\to[11221\underline{-22}1]\to [112211],
\end{gather*}
and a~cyclically adequate word with weight decreasing by~$4$.

In a~similar way one checks for the other seven cases in~\eqref{44} that one can apply Lemma~\ref{lm11}.
It is important to notice for later that in the f\/irst four cases, we do not need in fact that any of the neighboring
syllables of $\Delta^2$ is trivial.
\end{caselist}
\end{caselist}
\end{caselist}

\case There is no $\Delta^2$ subword in~$\beta$, but there are two~$\Delta$ subwords, i.e., we can apply Lemma~\ref{lm2}.
Thus
\begin{gather}
\label{ja}
\beta =\alpha_1\Delta\alpha_2\Delta
\end{gather}
with $\bar\alpha_1\alpha_2$ a~$\psi$-everywhere equivalent word.
If $\bar\alpha_1\alpha_2$ is central, the situation is clear.

\begin{caselist}
\case $\bar\alpha_1\alpha_2$ is split.
These are the words~$\beta$ (up to symmetries)
\begin{gather*}
\def\c#1{{\mathbf #1}}
\begin{array}{*3{@{}l@{
\qquad
}}c} \big[|1\c 21|2^k|121|1^l\big], & \big[|1\c 21|2^k|212|1^l\big], & \big[|212|2^k|1\c 21|1^l\big], & \big[|212|2^k|21\c 2|1^l\big].
\end{array}
\end{gather*}
They are distinguished by smoothings of the indicated boldfaced letter (giving as in Case~\ref{c1q} a~$(2,n)$-torus
link) and some letter outside the copies of~$\Delta$ (giving a~link of braid index~3).

\case We have
\begin{gather}
\label{aapp}
\alpha=\bar\alpha_1\alpha_2\doteq \big(\sigma_1^l\sigma_2^l\big)^k.
\end{gather}
We can obviously assume, by excluding $\Delta^2$ subwords, that none of $\alpha_i$ is the trivial (empty) word.

Recall the sliding~\eqref{sld} we used to bring two subwords~$\Delta$ together to form a~$\Delta^2$:
\begin{gather}
\label{sle}
\beta =\alpha_1\Delta\alpha_2\Delta\to\Delta\bar\alpha_1\alpha_2\Delta\doteq\Delta^2 \bar\alpha_1\alpha_2.
\end{gather}

\begin{caselist}
\case If one of the~$\Delta$ in~\eqref{ja} has neighbored letters of dif\/ferent index, after the sliding the
other~$\Delta$ close, will have neighbored letters to $\Delta^2$ of the same index.

Look at the words in the f\/irst row of~\eqref{44}.
These are built around the the four words for $\Delta^2$ factoring into two subwords of~$\Delta$.
In all four cases, the indicated to-change crossings lie entirely in one of the copies of~$\Delta$ in $\Delta^2$.
By symmetry, it can be chosen in either of them, in particular in the one we did not slide.

Thus one can undo in the same way bringing together the two copies of~$\Delta$ after either crossing changes (on the
right of~\eqref{sle}), and realize the two crossing changes in the original braid word~$\beta$ (on the left
of~\eqref{sle}).
The crossing changes in~\eqref{44} thus apply also in~$\beta$ to give positive words satisfying the assumptions of
Lemma~\ref{lm11}, and we are again done.

\case Now either of~$\Delta$ has neighboring syllables in $\alpha_i$ of the same index.

We assumed (by excluding $\Delta^2$ subwords in the present case) that none of $\alpha_i$ is a~trivial word.

\begin{caselist}
\case Let us f\/irst exclude that an $\alpha_i$ has length~$1$.
We will return to this situation later.

Since $[1|121|12]=[1|212|12]$ is central, none of the neighboring syllables can be trivial.
In particular, we can exclude the cases $l=1$, i.e., that $\alpha=\bar\alpha_1\alpha_2\doteq[12]^k$ in~\eqref{aapp}.

Sliding one copy of~$\Delta$ next to the other, we have a~word for $\Delta^2$ factoring into two words for~$\Delta$.
We can up to reversal assume that we slide the second~$\Delta$, while the f\/irst is f\/ixed.
As long as we switch crossings outside the second~$\Delta$ (i.e., in the f\/irst~$\Delta$, or outside~$\Delta^2$), we can
realize these switches also in~$\beta$ in~\eqref{ja} before the sliding.
We will now choose two such crossing changes so as to apply Lemma~\ref{lm11}.

Changing a~crossing in a~(non-trivial, as we assumed $l\ge 2$ in~\eqref{aapp}) syllable of $\alpha=\bar\alpha_1\alpha_2$
gives a~positive word $\beta '=\Delta^2\alpha'$ with a~trivial syllable in one of the~$\Delta$ and weight decreased by at most~2.

We compare such a~word with one obtained from the crossing changes indicated in Case~\ref{c122_}.
In all of $[121121]$, $[121212]$, $[212121]$ and $[212212]$ (the words for $\Delta^2$ factoring into two words
for~$\Delta$), one can change a~crossing among the f\/irst three letters to obtain $[1^22^2]$ or $[2^21^2]$, giving
a~cyclically adequate word $\tilde\beta $.
Up to $\sigma_1\leftrightarrow\sigma_2$ the above four options for $\Delta^2$ words reduce to two.

{\def\b#1{{\mathbf #1}}\def\i#1{{\mathit #1}} For $\Delta^2=[121212]$, we use the switch $[1\b 21212]$, giving
$[1^22^2]$, and a~cyclically adequate word~$\tilde\beta $ of weight by 4 less.
In $\beta '$ the weight decreases by at most 2, and thus Lemma~\ref{lm11} applies.

For $\Delta^2=[2\b 12212]$, we obtain $[2^21^2]$.
If it decreases the weight by 4, we are done.
Otherwise it decreases the weight by 2, and the letter following $\Delta^2$ is `2': $[|212212|2]$.

If the letter preceding $\Delta^2$ is also `2', i.e.~$[2|212212|\i 2]$, the italic letter switched decreases the weight
by at most 2, and leaves an isolated syllable `1' in $\Delta^2$.
Otherwise, the letter preceding~$\Delta^2$ is~`1'.
By remembering~\eqref{aapp} the word surrounding $\Delta^2$ and our assumption $l\ge 2$, we see that we must have
$[11|21\i 2212|2]$.
Then the italic letter switched decreases the weight by 2 and leaves an isolated~`2'.}

\case We must treat extra the case that one of $\alpha_i$ has one letter.

If both $\alpha_i$ have, these are the words
\begin{gather*}
[1|121|1|121],
\qquad
[1|121|1|212],
\qquad
[1|212|1|212].
\end{gather*}
These are knot diagrams of 8 crossings, and it is easily checked (as done in~\cite{evrdiff}) that only the third
is~$\psi$-everywhere equivalent, as desired.
(Here we use that for $\le 8$ crossings positive braid knots have all dif\/ferent Jones polynomial.)

When, say, $\alpha_2$ has more than one letter, these are up to interchanges, and (cyclically) permuting the last letter
of~$\alpha$ to the left, the following situations:
\begin{gather*}
[1|121|1|121|1],
\qquad
[1|121|1|212|1],
\qquad
[1|212|1|212|1].
\end{gather*}
In order removing the two copies of~$\Delta$ to give a~word~$\alpha$ as in~\eqref{aapp} up to cyclic permutations, we
need that $l=1$ and $k\nmid 3$.
But we excluded these cases already with the argument that then~$\beta$ has a~$\Delta^2$ subword.
\end{caselist}
\end{caselist}
\end{caselist}

\case If we cannot apply Lemma~\ref{lm1} or~\ref{lm2}, then we write up to cyclic permutations
\begin{gather}
\label{P1}
\beta =\Delta\alpha,
\end{gather}
and we assume that~$\alpha$ has no~$\Delta$ subword.
Let us assume also that~$\alpha$ has at least three letters.
(The other low crossing cases are rather obvious.)

We know, by case assumption, that a~trivial syllable in~$\alpha$ in~\eqref{P1} can occur only with its f\/irst or last letter.
We argue that, up to excluded cases, one can write~$\beta$ as in~\eqref{P1} so that even these two syllables are
non-trivial, that is, so that~$\alpha$ is adequate.

If we can ascertain such~$\alpha$, we are done with Lemma~\ref{lm11} as follows.
We compare two smoothings of~$\beta$.
On the one hand, smoothing any crossing in~$\alpha$ will give a~positive word of equal weight.
On the other hand, smoothing the non-repeating letter in $[121]$ or $[212]$ representing~$\Delta$ in~\eqref{P1} will
give a~cyclically adequate word of smaller weight.

We argue thus how to f\/ind a~presentation~\eqref{P1} of~$\beta$ with adequate (and not only cyclically adequate)~$\alpha$.

By permuting the last letter of~$\alpha$ to the left, we have some of the following situations (up to reversal and
$\sigma_1\leftrightarrow\sigma_2$):
\begin{gather*}
[1|121|1],
\qquad
[1|121|2],
\qquad
[2|121|2].
\end{gather*}
In either the f\/irst and third case, none of the neighboring syllables to~$\Delta$ can be trivial, because otherwise we
have a~$\Delta^2$ subword up to cyclic permutations, and we dealt with this case.
Thus permuting back all of~$\alpha$'s letters to the right, we have the shape we wanted.

In the second case, if the neighboring syllable `$2$' is trivial, then we have $[1|121|21]$.
Now choosing a~better~$\Delta$: $[11|212|1]$, one arrives (up to $\sigma_1\leftrightarrow\sigma_2$) in the third case above.
Thus if we cannot obtain the desired shape of~$\alpha$, the syllable after~$\Delta$ must be non-trivial: $[1|121|22]$.
Then compare the smoothings of the f\/irst and last `$2$': $[1|11|22]$ and $[1|121|2]$.
The f\/irst subword gives a~cyclically adequate word of smaller weight, the second one a~word of equal weight.
\end{caselist}

This f\/inishes the proof of Theorem~\ref{th1}.

\begin{rem}
The use of the Jones polynomial means \em{a priori} that we distinguish the links of $D_i'$ as links \em{with}
orientation, up to simultaneously reversing orientation of all components.
However, this restriction is not necessary, and in fact, we may see the links of $D_i'$ non-isotopic as unoriented
links.

Namely, in the non-positive case (family~\ref{type1} in Theorem~\ref{th1}), we have only knot diagrams~$D$, where the
issue is irrelevant.
In the positive case (families~\ref{type2} and~\ref{type3}), one can easily see that all $\beta '$ can be simplif\/ied to
a~positive braid word of two fewer crossings.
It is a~consequence of the minimal degree of the Jones polynomial (see, e.g.,~\cite{Fiedler}) and its reversing property
(see, e.g.,~\cite{LickMil}), that if the closures of two positive braids of the same exponent sum (and same number of
strings) are isotopic as unoriented links, then they are isotopic (up to reversing all components simultaneously) with
their positive orientations.
\end{rem}

\begin{rem}
Note also that, for links, our method does not restrict us to (excluding) isotopies between the links of $D_i'$ mapping
components in prescribed ways.
For example, the diagram~$D$ gives a~natural bijection between components of $D_i'$ and $D_j'$, but this correspondence
never played any role.
\end{rem}

\subsection{On everywhere dif\/ferent diagrams}

As an epilogue, we make a~useful remark that by modifying Lemmas~\ref{lm1} and~\ref{lm2}, one can easily see that the
construction of \em{everywhere different} diagrams (where all crossing-switched versions represent dif\/ferent links, and
here even in the strict, oriented, sense) is meaningless for 3-braids.

\begin{propos}
\label{p2}
There exist no everywhere different $3$-braids $($where we consider links up to oriented isotopy$)$.
\end{propos}

\begin{proof}
We show that for each 3-braid word~$\beta$, there are two crossing-switched versions giving conjugate braids.

We show this by induction over the word length.
Assume~$\beta$ is an everywhere dif\/ferent 3-braid word.

Obviously, by induction, we can restrict ourselves to braid words~$\beta$ with no $\sigma_i^{\pm 1}\sigma_i^{\mp 1}$
(whose deletion preserves the everywhere dif\/ferent property).

Thus~$\beta$ has an exponent vector, and for evident reasons, all syllables must be trivial.
(In particular the word length is even, and the closure is not a~2-component link.) If~$\beta$ were an alternating word,
then $\beta =(\sigma_1\sigma_2^{-1})^k$, which is obviously not everywhere dif\/ferent (at most two dif\/ferent links occur
after a~crossing change).

Since~$\beta$ is thus not alternating, it must contain cyclically a~word for $\Delta=[121]$ or~$[212]$, or~$\Delta^{-1}$.
Moreover, one easily sees that in a~subword $[1212]$ the edge crossing changes give the same braid, so that both
syllables around a~$\Delta$ (resp.\ $\Delta^{-1}$) must be negative (resp.\ positive).
In particular, dif\/ferent subwords of~$\beta$ representing~$\Delta$ (or $\Delta^{-1}$) are disjoint.

Again by symmetry reasons, there must be more than one $\Delta^{\pm 1}$ word.
(The braids $[121(-21)^k$ $-2]$ are not a~problem to exclude: switch the two crossings cyclically neighbored to $[121]$; if
$k=0$ there are two unknot diagrams.) Thus
\begin{gather*}
\beta =\Delta^{\pm 1}\alpha_1\Delta^{\pm 1}\alpha_2.
\end{gather*}
The word $\bar\alpha_1\alpha_2$ has now by induction two crossing changes (of the same sign) giving conjugate braids.
Since $\Delta^{\pm 1} \Delta^{\pm 1}$ is central (if not trivial), these crossing changes will remain valid in
$\Delta^{\pm 1}\Delta^{\pm 1}\bar\alpha_1\alpha_2$, and then also, by the sliding argument, in~$\beta$.
\end{proof}

\subsection*{Acknowledgment}

I wish to thank K.~Taniyama and R.~Shinjo for proposing the problems to me, and the referees for their helpful
comments.

\pdfbookmark[1]{References}{ref}
\LastPageEnding

\end{document}